\documentclass{amsart}
\usepackage{amsfonts}

\begin{document}
%\usepackage[inactive]{srcltx} % SRC Specials for DVI Searching

% Over-full v-boxes on even pages are due to the \v{c} in author's name
\vfuzz2pt % Don't report over-full v-boxes if over-edge is small

% THEOREM Environments ---------------------------------------------------
\newtheorem{thm}{Theorem}[section]
\newtheorem{cor}[thm]{Corollary}
\newtheorem{lem}[thm]{Lemma}
\newtheorem{prop}[thm]{Proposition}
\theoremstyle{definition}
\newtheorem{defn}[thm]{Definition}
\theoremstyle{remark}
\newtheorem{rem}[thm]{Remark}
\numberwithin{equation}{section}
% \newtheorem{thm}{Theorem}[subsection]
%\newtheorem{cor}[thm]{Corollary}
% \newtheorem{lem}[subsection]{Lemma}
 %\newtheorem{prop}[subsection]{Proposition}
 %\theoremstyle{definition}
 %\newtheorem{defn}[subsection]{Definition}
 %\theoremstyle{remark}
 %\newtheorem{rem}[subsection]{Remark}
 %\numberwithin{equation}{subsection}
% MATH -------------------------------------------------------------------
\newcommand{\CC}{\mathbb{C}}
\newcommand{\KK}{\mathbb{K}}
\newcommand{\ZZ}{\mathbb{Z}}
\def\a{{\alpha}}

\def\b{{\beta}}

\def\d{{\delta}}

\def\g{{\gamma}}

\def\l{{\lambda}}

\def\gg{{\mathfrak g}}
\def\cal{\mathcal }

\title{ On Nilpotent Leibniz Superalgebras}

\author{J.R. G{\'o}mez}
\author{R.M. Navarro}
\author{B. A. Omirov}

\address{Jos{\'e} Ram{\'o}n G{\'o}mez.\newline \indent
Dpto. Matem{\'a}tica Aplicada I, Universidad de Sevilla, Sevilla
(Spain)}

\email{jrgomez@us.es}

\address{Rosa Mar{\'\i}a Navarro.\newline \indent
Dpto. de Matem{\'a}ticas, Universidad de Extremadura, C{\'a}ceres
(Spain)}

\email{rnavarro@unex.es}

\address{B. A. Omirov.\newline \indent
Institute of Mathematics (Uzbekistan)}

\email{omirovb@mail.ru}

%\subjclass{17B30}

%\keywords{Leibniz superalgebras, nilindex}

%\begin{abstract}

 %\end{abstract}
 \thanks{This work was supported by the PAICYT,
 FQM143 of the Junta de Andalucia (Spain), by INTAS (Ref. Nr. 04-83-3035) and by Junta de
Extremadura-Consejer\'{\i}a de Infraestructuras y Desarrollo
Tecnol\'ogico and Feder (Ref. Nr. 3PR05A074)}

%%% ----------------------------------------------------------------------
\maketitle
%%% ----------------------------------------------------------------------

\begin{abstract}

The aim of this work is to present the first problems that appear
in the study of nilpotent Leibniz superalgebras. These
superalgebras and so the problems, will be considered as a natural
generalization of nilpotent Leibniz algebras and Lie
superalgebras.

 \end{abstract}

%%% ----------------------------------------------------------------------
\maketitle
%%% ----------------------------------------------------------------------


\section{Introduction}

The notion of Leibniz superalgebras was firstly introduced in
\cite{AO0}, although graded Leibniz algebra was considered before
in work \cite{Liv}. As Leibniz algebras are a generalization of
Lie algebras \cite{Loday}, then  many of the features of Leibniz
superalgebras are generalization of Lie superalgebras.

The study of nilpotent Leibniz algebras \cite{AO0}, \cite{AO1},
\cite{AO2} shows that many nilpotent properties of Lie algebras
can be extended for nilpotent Leibniz algebras. The results of
nilpotent Leibniz algebras may help us to study nilpotent Leibniz
superalgebras. However, nilpotent Leibniz superalgebras turn out
more complex than nilpotent Lie superalgebras.

This is the frame of our work: nilpotent Leibniz superalgebras. In
the case of Leibniz algebras appears the notion of zero-filiform
algebra, this notion does not exist in Lie algebras. This algebra
has maximal nilindex. In Leibniz superalgebras we offer the
analogue of zero-filiform Leibniz algebras, that is zero-filiform
Leibniz superalgebras, ${\cal ZF}^{n,m}$. But not all of
zero-filiform Leibniz superalgebras have maximal nilindex, there
is only one in particular pair of dimensions: $(n,n)$, $(n,n+1)$,
as could be seen in theorem \ref{ayupov}.

Before to studying general classes of Leibniz superalgebras
(zero-filiform and filiform Leibniz superalgebras) we had to solve
the problem of finding a suitable basis; a so-called adapted
basis, see theorems \ref{adapted1} and \ref{adapted2}.
 The function $f(n,m)$, defined as the maximal nilindex for Lie
superalgebras of type (n,m), is always $\leq$ $n+m-1$ \cite{GKN}
but the same function $f(n,m)$ on Leibniz superalgebras can be
$n+m$. Also, this value is only obtained for the special
zero-filiform Leibniz superalgebras that appear in the particular
pair of dimensions mentioned above. By applying direct sum with
$\CC$ to these special Leibniz superalgebras we obtain that
$f(n+1,n)=2n$ and $f(n,n+2)=2n+1$, see theorem \ref{index1}.
Maximal nilindex determination is an open problem for the general
case.

Analogously as for Lie superalgebras \cite{GKN} we will refer the
nilpotent Leibniz superalgebras of type $(n,m)$, with nilindex
$f(n,m)$, as maximal class Leibniz superalgebras. We will denote
the variety of these Leibniz superalgebras as ${\cal M}^{n,m}$.

In this work we have obtained many  results concerning nilindex
and the function maximal nilindex $f(n,m)$; we have found many of
superalgebras with open orbits which determine irreducible
components of the variety of nilpotent Leibniz superalgebras, and
we have obtained the relative position of the subvarieties ${\cal
M}^{n,m}$ and ${\cal ZF}^{n,m}$ in some cases. Also we conjecture
that there exists one unique Leibniz superalgebra of type
$(n+1,n)$ with nilindex equal to $2n$ (conjecture 2).

In this paper, most of classification proofs are omitted because
of they are very laborious and they do not apport any new idea.

\section{Preliminaries}

The vector space V is said to be $\ZZ_2$-graded if it admits a
decomposition in direct sum, $V=V_0\bigoplus V_1$. An element X of
V is called homogeneous of degree $\gamma$, $\gamma\in \ZZ_2$, if
it is an element of $V_{\gamma}$ . In particular, the elements of
$V_0$ (resp. $V_1$) are also called even (resp. odd).

Let $V=V_0\bigoplus V_1$ and $W=W_0\bigoplus W_1$ be two
$\ZZ_2$-graded vector spaces. A linear mapping $f : V\rightarrow
W$ is said to be homogeneous of degree $\gamma$, $\gamma\in
\ZZ_2$, if $f(V_{\alpha})\subseteq W_{\alpha+\gamma(mod 2)}$ for
all $\alpha\in \ZZ_2$. In particular, if the linear mapping is
homogeneous of degree $0$ is said to be a  homomorphism of the two
$\ZZ_2$-graded vector spaces. Now it is clear how we define an
isomorphism or an automorphism of $\ZZ_2$-graded vector spaces.

We say that two Leibniz superalgebras, $L_1$ and $L_2$, are
isomorphic if there exists a $\ZZ_2$-graded vector space
isomorphism, $\varphi : L_1\rightarrow L_2$, satisfying
$\varphi([X,Y])=[\varphi(X),\varphi(Y)]$ for all X, Y of $L_1$.
Then $\varphi$ is called an isomorphism of Leibniz superalgebras
and it is always assumed to be consistent with
$\ZZ_2$-graduations; that is, they are homogeneous linear mappings
of degree zero.

\begin{defn}\label{defB} \cite{AO0}. A $\ZZ_2$-graded vector space $L=L_0\bigoplus L_1$ is
called a Leibniz superalgebra if it is equipped with a product
$[-,-]$ which satisfies the following conditions:
$$[L_{\alpha},L_{\beta}]\subseteq L_{\alpha+\beta(mod2)} \ for \
all \ \alpha,\beta\in \ZZ_2$$
$$[x,[y,z]]=[[x,y],z]-(-1)^{\alpha\beta}[[x,z],y] -
 \ graded \ Leibniz  \ identity$$ for all $x\in L$, $y\in L_{\alpha}$,
$z\in L_{\beta}$, $\alpha,\beta\in \ZZ_2$.
\end{defn}

Note that if a Leibniz superalgebra L satisfies the identity
$[x,y]=-(-1)^{\alpha\beta}[y,x]$ for all $x\in L_{\alpha}$, $y\in
L_{\beta}$ , then the graded Leibniz identity  is reduced to the
following graded Lie identity:
$$(-1)^{\alpha\gamma}[x,[y,z]]+
(-1)^{\alpha\beta}[y,[z,x]]+(-1)^{\beta\gamma}[z,[x,y]]=0.$$

Therefore Leibniz superalgebras are a generalization of Lie
superalgebras.

If we denote by $R_X$ the right multiplication operator, i.e. $R_X
: L\rightarrow L$, then the graded Leibniz identity can be
expressed in the following form:
$$R_{[X,Y]}=R_{Y}R_{X}-(-1)^{\alpha\beta}R_XR_Y \ \ \ \ \ (1)$$
where $X\in L_{\alpha}, Y \in L_{\beta}$.

We denote by $R(L)$ the set of all right multiplication operators.
It is not difficult to prove that $R(L)$ with the multiplication
defined by:
$$<R_a, R_b>:=R_aR_b-(-1)^{\alpha\beta}R_bR_a \ \ \ \ \ (2)$$
for $R_a\in {R(L)_{\alpha}}$, $R_b\in {R(L)_{\beta}}$, becomes a
Lie superalgebra.

In order to provide an example of non-Lie Leibniz superalgebra, we
can consider an associative superalgebra, $A=A_0\bigoplus A_1$,
and a linear mapping $D : A\rightarrow A$ satisfying the
condition:
$$D(a(Db))=DaDb=D((Da)b)$$ for all $a, b\in A$
and define a new multiplication over the underlying $\ZZ_2$-graded
vector space, $< , >$, by:
$$<a,b>_D:=a(Db)-(-1)^{\alpha\beta}D(b)a$$ for $a\in A_{\alpha}, b\in
A_{\beta}$ . Then $A$ equipped with multiplication $< , >$ becomes
a Leibniz superalgebra, which in general is not a Lie
superalgebra.

The descending central sequence of a Leibniz superalgebra
$L=L_0\bigoplus L_1$ is defined by $C^0(L)=L$,
$C^{k+1}(L)=[C^k(L),L]$ for all $k\geq 0$. If $C^k(L)={0}$ for
some $k$, the Leibniz superalgebra is called nilpotent.  The
smallest natural number $k$ such as $C^k(L)={0}$ is called the
nilindex of $L$.

We denote by $N^{n,m}$ the variety of nilpotent Leibniz
superalgebras $L=L_0\bigoplus L_1$ with $dimL_0=n$, $dimL_1=m$;
and by $Leib^{n,m}$ denote the variety of Leibniz superalgebras.
The above property $C^k(L)={0}$ can be realized via finite numbers
of polynomial relations on structure constants and therefore, this
set forms a subvariety of the variety $Leib^{n,m}$.

For Leibniz superalgebras we have the following analogue of
Engel's theorem:

\begin{thm}  {\bf (Engel's theorem)} \cite{AO0}. A finite dimensional Leibniz
superalgebra $L$ is nilpotent if and only if the operators $R_X$
are nilpotent for all $X\in L$.
\end{thm}

If we take an homogeneous basis $\{X_{0}, ... , X_{n-1}, Y_{1},
.... , Y_{m}\}$ for $L$ $(L\in Leib^{n,m})$, the superalgebra is
completely determined by:

$$[X_i,X_j]=\sum\limits_{k=1}^{n-1}C_{ij}^{k}X_{k}, \ \ \ \
[X_i,Y_j]=\sum_{k=1}^{m}D_{ij}^{k}Y_{k},$$
$$[Y_i,X_j]=\sum_{k=1}^{m}E_{ij}^{k}Y_{k}, \ \ \ \
[Y_i,Y_j]=\sum_{k=1}^{n-1}F_{ij}^{k}X_{k},$$ where $\{C_{ij}^{k},
D_{ij}^{k}, E_{ij}^{k}, F_{ij}^{k}\}$ are structure constants.
These structure constants verify the restrictions obtained by the
graded Leibniz identity \cite{AO0}.

Let $V=V_0\bigoplus V_1$ be the underlying vector space of $L$,
$L=L_0\bigoplus L_1\in Leib^{n,m}$ and let $G(V)$ be the group of
the invertible linear mappings of the form $f=f_0+f_1$ such that
$f_0\in GL(n,C)$ and $f_1\in GL(m,C)$ $(G(V)=GL(n,C)\bigoplus
GL(m,C))$. The action of $G(V)$ on $Leib^{n,m}$ induces an action
on the Leibniz superalgebras variety: two laws $\lambda_1 \  and \
\lambda_2$ are isomorphic, if there exists a linear mapping $f$,
$f=f_0+f_1\in G(V)$, such that
$$\lambda_2(X,Y)=f_{\alpha+\beta}^{-1}(\lambda_1(f_{\alpha}(X),
f_{\beta}(Y))) \ \ \ for \ \ all \ \ X\in V_{\alpha}, Y\in
V_{\beta}.$$

We denote by $O(\lambda)$ the orbit of $\lambda$ corresponding to
this action.

The superalgebras with open orbits in $N^{n,m}$ are called rigid.
The closures of these open orbits give irreducible components of
the variety $N^{n,m}$. Then, the fact of finding such algebras is
crucial for the description of the variety $N^{n,m}$.

The description of the variety of any class of algebras or
superalgebras is a difficult problem. Different papers, for
example \cite{AOR}, \cite{BS}, \cite{GO}, \cite{GK} are concerning
the applications of algebraic groups theory to the description of
the variety of Lie algebras.

\begin{defn} For a Leibniz superalgebra $L=L_0\bigoplus L_1$ we define the
set $Z(L)$, $Z(L)=\{X\in L  \ : \ [L,X]=0\}$ which  will be called
the right annihilator of $L$.
\end{defn}

It is easy to see that $Z(L)$ is a two-sided ideal of $L$ and
$[X,X]\in Z(L)$ for any $X\in L_0$, this notion is good and
compatible with the right annihilator in Leibniz algebras. If we
consider $I=ideal<[X,Y]+(-1)^{\alpha\beta}[Y,X] \ : \ X\in
L_{\alpha}, Y\in L_{\beta}>$, then $I\subseteq Z(L).$

\begin{defn}
For a Leibniz superalgebra $L$ we define the sets
$$L(L)=\{X\in L   :  [X,L]=0\},$$
$$Cent(L)=\{X\in L   :  [X,L]=[L,X]=0\}$$
which are called left annihilator and center of $L$, respectively.
\end{defn}

We can extract a result of Leibniz algebras \cite{AOR} and apply
it to Leibniz superalgebras : for any $s, r\in \mathbb{N}$ the
following subsets of $Leib^{n,m}$ are closed relatively to the
Zariski topology:
\begin{enumerate}
\item $ \{\mu\in Leib^{n,m}\enskip |\enskip dim\mu^s\le r\}$

\item $ \{\mu\in Leib^{n,m}\enskip |\enskip dimZ(\mu)\ge s\}$

\item $ \{\mu\in Leib^{n,m}\enskip |\enskip dimL(\mu)\ge s\}$

\item $ \{\mu\in Leib^{n,m}\enskip |\enskip dimCent(\mu)\ge s\}$

\end{enumerate}

where $\mu^s=C^{s}(\mu)$. Hence, a superalgebra $\mu$ does not
belong to $clO(\lambda)$ if one of the following conditions holds:
\begin{enumerate}
\item  $dim\lambda^s<dim\mu^s$ for some $s,$ \item $
dimZ(\lambda)>dimZ(\mu),$ \item $dimL(\lambda)>dimL(\mu),$ \item $
dimCent(\lambda)>dimCent(\mu).$
\end{enumerate}

\section{Details in Leibniz superalgebras}

Let $L=L_0\bigoplus L_1$ be a nilpotent Leibniz superalgebra with
$dimL_0=n$ and $dimL_1=m.$ From $(2)$ we have that $R(L)$ is a Lie
superalgebra, in particular $R(L_0)$ is a Lie algebra. As $L_1$
has $L_0$-module structure we can consider $R(L_0)$ as a subset of
$GL(V)$, where $V$ is vector space corresponding to $L_1$. So, we
have a Lie algebra formed by nilpotent endomorphisms of $V$. And
applying the Engel's theorem \cite{Jac} we have existence of
subspaces of $V$:
$$V_0\subseteq V_1\subseteq V_2\subseteq ... \subseteq V_m=V,$$ with $dim(V_i)=i$
and $R(L_0)(V_{i+1})\subset V_i.$

We define two new descending sequences, $C^k(L_0)$, $C^k(L_1)$ as
follows: $C^0(L_i)=L_i, \ \ C^{k+1}(L_i)=[C^k(L_i),L_0] \ for \
k\geq 0, \ i\in \{0,1\}.$ Analogously as for Lie superalgebras
\cite{GKN}, if $L=L_0\bigoplus L_1$ is a nilpotent Leibniz
superalgebra, then $L$ has super-nilindex or s-nilindex $(p,q)$ if
the following conditions hold:
$$(C^{p-1}(L_0))(C^{q-1}(L_1)) \neq 0, \qquad C^p(L_0)=C^q(L_1)=0$$

We have for Lie superalgebras the invariant called characteristic
sequence that can be naturally extended  for Leibniz
superalgebras. Thus, we have the following definition.

\begin{defn} For an arbitrary element $X\in L_0$, the operator $R_X$ is
a nilpotent endomorphism of space $L_i$, where $i\in \{0,1\}$. We
denote by $gz_i(X)$ descending sequence of dimensions of Jordan
blocks of $R_X$. We define the invariant of a Leibniz superalgebra
$L$ as follows:
$$
gz(L)=\left(\left.\max_{X\in L_0-[L_0,L_0]} gz_0(X) \ \right|
\max_{\widetilde{X}\in L_0-[L_0,L_0]} gz_1(\widetilde{X})\right),
$$
where $gz_i$ is the lexicographic order.

The couple $gz(L)$ is called characteristic sequence of Leibniz
superalgebra $L$.
\end{defn}

We denote by $N^{n,m}_{p,q}$  the subset of the set $N^{n,m}$, with
s-nilindex $(k_0,k_1)$, where $k_0\leq p$ and $k_1\leq q$.

\begin{lem} The set $N^{n,m}_{p,q}$ is a subvariety of the variety $N^{n,m}.$
\end{lem}
\begin{proof} The proof of this lemma is evident, because the set $N^{n,m}_{p,q}$
can be realized via finite numbers of polynomial equations of the
structure constants.
\end{proof}

\begin{defn} A Leibniz superalgebra $L$,  $ L \in N^{n,m}$, is  called zero-filiform
if its s-nilindex is $(n,m)$.
\end{defn}

We denote by ${\cal ZF}^{n,m}$  the set of zero-filiform Leibniz
superalgebras.

\begin{rem}

\

\begin{itemize}
\item[$1)$] If $L=L_0\bigoplus L_1$  is a zero-filiform Leibniz
superalgebra then from \cite{AO1} we have that $L_0$ is a
zero-filiform Leibniz algebra.

\item[$2)$] Since ${\cal ZF}^{n,m}=N^{n,m}\backslash(N_{n-1,m}^{n,m}\cup
N_{n,m-1}^{n,m})$, then ${\cal ZF}^{n,m}$ is an open set in
Zariski topology.

\item[$3)$] Note that zero-filiform Leibniz superalgebra $L$ of
type $(n,m)$ can be realized as superalgebra with $gz(L)=(n|m).$
\end{itemize}
\end{rem}

Before to studying general classes of Leibniz (super)algebras, it
is useful to solve the problem of finding a suitable basis; a
so-called adapted basis. This question is not trivial even for Lie
superalgebras and it is difficult to demonstrate the general
existence of such a basis for Leibniz superalgebras. Particularly,
we prove in the following theorem that there always exists an
adapted basis for the class of zero-filiform Leibniz
superalgebras.

\

\begin{thm}\label{adapted1} If $L=L_0\bigoplus L_1 \in {\cal ZF}^{n,m}$, then there
exists an adapted basis of $L$, namely $\{X_0,X_1, ... , X_{n-1},
Y_1, Y_2, ..., Y_m\}$, with $\{X_0,X_1, .... , X_{n-1}\}$ a basis
of $L_0$ and $\{Y_1, Y_2, ..., Y_m\}$ a basis of $L_1$, such that:
$$\begin{array}{llc}
[X_i,X_0]=X_{i+1},& 0 \leq i \leq n-2, & [X_{n-1},X_0]=0,\\[1mm]
[Y_j,X_0]=Y_{j+1},& 1 \leq j \leq m-1,& [Y_m,X_0]=0. \\[1mm]
\end{array}$$
Moreover $[Y_j,X_k]=0$ for $1 \leq j \leq m$ and  $1 \leq k \leq
n-1$, and the omitted products of $L_0=<X_0,X_1, ... , X_{n-1}>$
vanish.
\end{thm}

\begin{proof} As $L=L_0\bigoplus L_1$ is a zero-filiform Leibniz superalgebra,
then $L_0$ is a zero-filiform Leibniz algebra. Thus from
\cite{AO1} we
have an adapted basis for $L_0:$ $\{X_0, \ X_1, \ ..., \\
X_{n-1}\}$ with $[X_i,X_0]=X_{i+1},$  $0\leq i\leq {n-2},$
$[X_{n-1},X_0]=0, \ [X_i,X_k]=0$ for $0\leq i\leq {n-1},$  and
$1\leq k\leq n-1$.

As we stated at the beginning of this section we have:
$$0\subset V_1\subset ... \subset V_m \ \mbox{ with } \ dim(V_{i+1}/V_i)=1,$$
where each $V_i$ is the vector space of generators: $\{Y_1, Y_2,
...., Y_i\}$, $V_i=<Y_1, Y_2, ..., Y_i>$, with $[V_{i+1},L_0]=V_i.$

So, $[V_1,L_0]=0$ and then $[Y_1,X_i]=0 \ \forall i$. As
$[V_2,L_0]=V_1$ we have that exists a non-null scalar namely
$\l_{2}$ such that $[Y_2,X_{i_2}]=\l_{2} Y_1$. By induction, it is
possible to prove that there exists a set of non-null scalars
$\{\l_{2},\l_{3}, \dots, \l_{m}\}$ and vectors
$\{X_{i_2},X_{i_3},\dots,X_{i_m}\}\subseteq
\{X_{0},X_{1},\dots,X_{n-1}\}$ verify
$$
[Y_k,X_{i_k}]=\l_{k}Y_{k-1}+\Psi_k(Y_{k-2},\dots,Y_1), \quad
 2 \leq k \leq m, \qquad (3)$$
where $\Psi_k(v_1,v_2,\dots,v_s)$ represents a linear combination
of the vectors $\{v_1,v_2,\dots,v_s \}$.

Using the graded Leibniz identity we can assert that $i_2=i_3=
\cdots =i_m=0$. In fact, if there exists $i_k \in \{1,
\dots,n-1\}$ we have
$$ [Y_k,X_{i_k}]=[Y_k,[X_{i_k-1},X_0]]=
\underbrace{ [ \underbrace{[Y_k,X_{i_k-1}]}_{\subseteq
V_{k-1}},X_0]}_{\subseteq V_{k-2}} -
\underbrace{[\underbrace{[Y_k,X_0]}_{\subseteq V_{k-1} } ,
X_{i_k-1}]}_{\subseteq V_{k-2}} $$ But from (3) we obtain that
$Y_{k-1} \in V_{k-2}$ which is a contradiction with the definition
of the subspaces $V_i$.

Thus, we have the following expression for the basis vectors
$$\begin{array}{ll}
[Y_1,X_j]=0, & \forall j
\\[1mm]
[Y_2,X_0]=\l_2 Y_1,& \l_2\neq 0
\\[1mm]
[Y_3,X_0]=\l_3 Y_2 + \Psi_3(Y_1),& \l_3 \neq 0
\\[1mm]
%[X_0,Y_4]=\l_4 Y_3 + \Psi_4(Y_2,Y_1)&\\ [1mm][X_1,Y_4]=\d_4 Y_3 +
%\Phi_4(Y_2,Y_1)& (\l_4,\d_4)\neq (0,0)
%\\[1mm]
%\hspace*{5em} \vdots& \\[1mm]
[Y_i,X_0]=\l_i Y_{i-1} + \Psi_i(Y_{i-2},\dots,Y_1),& \l_i \neq 0,
\ 4 \leq i \leq m
\end{array}
$$

Using the change of basis
$$\left\{ \begin{array}{ll}
X'_i=X_i, & 0 \leq i \leq n-1\\
Y'_m=Y_m&\\
Y'_{m-j}=[Y'_{m-j+1},X'_0], & 1 \leq j \leq m
\end{array}\right.$$
 and namely $Y_j=Y'_{m-j+1}$ for $1 \leq j \leq m$ we obtain $[Y_j,X_0]=Y_{j+1}$
 with $1 \leq j \leq m-1$ and $[Y_m,X_0]=0$. Only rest to prove
$[Y_j,X_k]=0$ for $1 \leq j \leq m$ and  $1 \leq k \leq n-1$ for
to conclude the proof.

Using graded Leibniz identity for the vectors $(Y_j,X_0,X_0)$ with
$1 \leq j \leq m$, we obtain that $[Y_j,X_1]=0$. It is easily seen
 by induction in $k$ and using graded Leibniz identity for the
vectors $(Y_j,X_{k-1},X_0)$ that $[Y_j,X_k]=0$ for $1 \leq j \leq
m$ and  $1 \leq k \leq n-1$ which concludes the proof.
\end{proof}

\begin{defn} A Leibniz superalgebra of $N^{n,m}$  is called filiform
if its s-nilindex is (n-1,m).
\end{defn}

\begin{rem}

\

\begin{itemize}
\item[$1)$] If $L=L_0\bigoplus L_1$  is a filiform Leibniz superalgebra
then from \cite{AO1} we have that $L_0$ is a filiform Leibniz
algebra.
\item[$2)$] Note that a filiform Leibniz superalgebra $L$ of type
$(n,m)$ can be realized as a Leibniz superalgebra with
$gz(L)=(n-1,1| m).$
\end{itemize}
\end{rem}

We denote by $F^{n,m}$  the set of all filiform Leibniz
superalgebras.

The next theorem shows that in the class of filiform Leibniz
superalgebras it is also possible to assert about existence of
an adapted basis.

\

\begin{thm}\label{adapted2} Let $L=L_0\bigoplus L_1$ be a filiform Leibniz
superalgebra, $L \in F^{n,m}$. Then there exists a basis $\{X_0,
X_1, ..., X_{n-1}, Y_1, Y_2, ..., Y_m\}$ of L such that $L$ can be
expressed in one of the following of laws:

$$(I):\left\{ \begin{array}{l}[X_i,X_0]=X_{i+1}, \quad 1\leq i\leq
{n-2}\\ [1mm]
[Y_j,X_0]=Y_{j+1}, \quad 1\leq j\leq {m-1} \\[1mm]
[Y_m,X_0]=0  \\[1mm]
[Y_j,X_k]=0, \quad  1\leq j\leq m, \ 2\leq k\leq {n-1}\\ [1mm]
[X_0,X_0]=X_2 \\[1mm]
[X_0,X_1]=\alpha_3X_3+...+\alpha_{n-2}X_{n-2}+\theta X_{n-1} \\[1mm]
[X_i,X_1]=\alpha_3X_{i+2}+...+\alpha_{n-i}X_{n-1}, \quad  1\leq
i\leq {n-3}\end{array}\right.$$

$$(II):\left\{ \begin{array}{l}[X_i,X_0]=X_{i+1}, \quad 2\leq i\leq
{n-2}\\ [1mm]
[Y_j,X_0]=Y_{j+1}, \quad 1\leq j\leq {m-1} \\[1mm]
[Y_m,X_0]=0  \\[1mm]
[Y_j,X_k]=0, \quad  1\leq j\leq m, \ 2\leq k\leq {n-1}\\ [1mm]
[X_0,X_0]=X_2 \\[1mm]
[X_0,X_1]=\beta_3X_3+...+\beta_{n-1}X_{n-1} \\[1mm]
[X_1,X_1]=\gamma X_{n-1}\\[1mm]
[X_i,X_1]=\beta_3X_{i+2}+...+\beta_{n-i}X_{n-1}, \quad 2\leq i\leq
{n-3}\end{array}\right.$$

$$(III):\left\{ \begin{array}{ll}
[X_i,X_0]=-[X_0,X_i]=X_{i+1}&1\leq i\leq n-2\\{}
[X_i,X_j]=-[X_j,X_i]\in
lin<X_{i+j+1},X_{i+j+2},\ldots,X_{n}>&1\leq i, j\leq n-1\\{}
[Y_j,X_0]=Y_{j+1}&1\leq j\leq m-1
\end{array}\right.
$$
In (I) and (II) the omitted products of $L_0=<X_0,X_1, ... ,
X_{n-1}>$ vanish.
\end{thm}

\begin{proof} If $L_0$ is a  non-Lie filiform Leibniz algebra, then from
\cite{AO1} we have that there exists a basis of $L_0:$ $\{X_0, \
X_1, \ ..., X_{n-1}\}$ such that $L_0$ can be expressed as the
even products (i.e. the products $[X_i,X_j]$) of (I) or as the
even products of (II). Then applying a similar reasoning as in
theorem \ref{adapted1} we obtain the existence of vectors of
$L_1:$ \ $\{Y_1, ...., Y_m\}$ which satisfy the multiplications
(I) and (II) of the theorem.

If $L_0$ is a filiform Lie algebra, then using the even products
from \cite{GK} and using the mentioned reasoning for $\{Y_1, \dots
, Y_m\}$ way we obtain the family (III). \end{proof}

In \cite{AO0} we have the following result
\begin{thm}\cite{AO0} \label{ayupov} Let L be a n-dimensional Leibniz superalgebra with maximal index of
nilpotency. Then L is isomorphic to one of the two following non
isomorphic algebras:
$$[e_i,e_1]=e_{i+1} \ 1\leq i\leq n-1,$$

$$[e_i,e_1]=e_{i+1}, \ 1\leq i\leq {n-1}, \ [e_i,e_2]=2e_{i+2}, \ 1\leq i\leq {n-2},$$
where omitted products are zero.
\end{thm}

Seeing the process of the proof of this theorem we observe that
the Leibniz superalgebra $$[e_i,e_1]=e_{i+1} \ 1\leq i\leq n-1$$
is really a split Leibniz superalgebra, i.e. is a Leibniz algebra
with all the basis elements $e_i$ even ones. So, this case is not
interesting for our study. However, the Leibniz superalgebra
$$[e_i,e_1]=e_{i+1}, \ 1\leq i\leq {n-1}, \ [e_i,e_2]=2e_{i+2}, \
1\leq i\leq {n-2}$$ has even part all the $e_i$ with $i$ even, and
the odd part is constituted by all the $e_i$ with $i$ odd. That
is, $L=L_0 \bigoplus L_1=< e_2, e_4, \dots > \bigoplus <
e_1,e_3,\dots >$, thus there are two possibilities:
\begin{itemize}
 \item if $n$ is even, then we have $L=<e_2, e_4, \dots,e_n> \bigoplus <e_1,e_3,\dots,e_{n-1}>$ so,
$dim(L_0)=dim(L_1)$.\\
 \item if $n$ is odd, then we have $L=<e_2,e_4,\dots,e_{n-1}> \bigoplus <e_1,e_3,\dots,e_{n}>$ so,
$dim(L_1)=1+dim(L_0)$.
\end{itemize}
As for Lie superalgebras \cite{GKN}, the function that gives the
maximal nilindex for each pair of dimensions $n$ and $m$
(dimensions of the even and odd parts, respectively) will be noted
by $f(n,m)$. Thus we have the following important theorems for our
study.

\

\begin{thm}
 $f(n,m)$ is equal to $n+m$ if and only if $m=n$ or
$m=n+1.$
\end{thm}

\

By applying direct sum with $\CC$ to the non-split superalgebras
of the theorem \ref{ayupov}, we obtain the following result.

\

\begin{thm}\label{index1} $f(n+1,n)=2n$ and $f(n,n+2)=2n+1$.
\end{thm}

\

At this point, for the rest of possibilities of the pair $(n,m)$
the function $f(n,m)$ is unknown.

Since the non-split superalgebras of the theorem \ref{ayupov} are
zero-filiform Leibniz superalgebras, in the following sections we
will start considering the set ${\cal ZF}^{n+1,m}.$ Moreover, in
the next sections we will study Leibniz superalgebras with the
dimension of the odd part up to three and generic dimension of the
even part, and generic dimension of the odd part and
two-dimensional even part.

\section{Leibniz Superalgebras with two-dimensional odd part}
We will consider the case $n=1$ separately.

\subsection{Case $n=1$}
\begin{lem}Let $L$ be
any Leibniz superalgebra $L \in {\cal ZF}^{2,2}$. Then it is
isomorphic to one of the following Leibniz superalgebras, pairwise
non-isomorphic, that can be expressed in an adapted basis
$\{X_0,X_1,Y_1,Y_2\}$ by
$$\begin{array}{ll}%\mu_{1} =\left\{\begin{array}{ll}
%[X_0,X_0]=X_{1},& \\[1mm]
%[Y_1,X_0]=Y_2 & \\[1mm]
%[Y_1,Y_1]=X_1 &
%\end{array}\right.&
\mu_{1}^{\a} =\left\{\begin{array}{ll}
[X_0,X_0]=X_{1},&  \\[1mm]
[Y_1,X_0]=Y_2 & \\[1mm]
[X_0,Y_1]=\a Y_2,& \a \in \CC\\[1mm]
[Y_1,Y_1]=X_1 &
\end{array}\right.&
\mu_{2} =\left\{\begin{array}{ll}
[X_0,X_0]=X_1& \\[1mm]
[Y_1,X_0]=Y_2 & \\[1mm]
[X_0,Y_1]=\frac{1}{2} Y_2& \\[1mm]
[Y_1,Y_1]=X_0 & \\[1mm]
 [Y_2,Y_1]=X_1 &
\end{array}\right.\\ \\

\end{array}
$$
\end{lem}

\begin{proof} By using a generic change of basis, along with the
graded Leibniz identity we obtain the lemma. \end{proof}

\begin{rem}

By using the change of basis $e_1=Y_1, \ e_2=X_0, \
e_3=\frac{1}{2}Y_2, \ e_4=\frac{1}{2}X_1$, it is easy to see that
$\mu_2$ is a superalgebra of the theorem \ref{ayupov} for the case
$(2,2)$.
\end{rem}

\begin{prop} $f(2,2)=4(=n+m)$
\end{prop}

\begin{proof} Is a corollary of the above lemma.
\end{proof}

%\noindent By using an equivalent notation to the one of \cite{GKN}
%we obtain that

\begin{prop}$${\cal M}^{2,2}=O(\mu_2)$$
\end{prop}

\begin{proof} According to the above lemma it is sufficient to prove that
if $L \in N^{2,2}$ with nilindex 4, then $L \in {\cal ZF}^{2,2}$.
In fact, if $L \notin {\cal ZF}^{2,2}$ then we have two possible
cases:

Case 1. $L=L_0+L_1$, with $L_0$ abelian. In this case $L$ is a Lie
superalgebra and then it will have nilindex $< 4$, see \cite{GKN}.

Case 2. $L=L_0+L_1$, with $L_0$ zero-filiform Leibniz algebra. In
this case as $L \notin {\cal ZF}^{2,2}$ then $Y_2 \notin {\cal
C}^1(L)$ ($\{Y_1,Y_2\}$ basis of $L_1$) which leads to nilindex $<
4$. \end{proof}
\begin{prop}The orbit $O(\mu_2)$ is a Zariski open subset of
$N^{2,2}$
\end{prop}
\begin{proof} It is a consequence of $${\cal M}^{n+1,m}=N^{n+1,m}-N^{n+1,m}_{f(n+1,m)-1}$$ \end{proof}

\subsection{Case $n>1$. }

\begin{lem}Let $L$ be
any zero-filiform Leibniz superalgebra $L \in {\cal ZF}^{n+1,2}$,
with $n\geq 2$. Then it is isomorphic to one of the following
Leibniz superalgebras, pairwise non-isomorphic, that can be
expressed in an adapted basis $\{X_0,X_1,\dots,X_n, Y_1,Y_2\}$ by

$$\begin{array}{ll}\mu_{1}^{\a} =\left\{\begin{array}{ll}
[X_i,X_0]=X_{i+1},& 0\leq i \leq n-1 \\[1mm]
[Y_1,X_0]=Y_2 & \\[1mm]
[X_0,Y_1]=\a Y_2,& \a \in \CC\\[1mm]
[Y_1,Y_1]=X_n &
\end{array}\right.&
\mu_{2} =\left\{\begin{array}{ll}
[X_i,X_0]=X_{i+1},& 0\leq i \leq n-1\\[1mm]
[Y_1,X_0]=Y_2 & \\[1mm]
%[X_0,Y_1]=\frac{1}{2} Y_2& \\[1mm]
[Y_1,Y_1]=X_{n-1} & \\ [1mm] [Y_2,Y_1]=X_{n} &
\end{array}\right.\\ \\
\mu_{3} =\left\{\begin{array}{ll}
[X_i,X_0]=X_{i+1},& 0\leq i \leq n-1\\[1mm]
[Y_1,X_0]=Y_2 & \\[1mm]
[X_0,Y_1]=-Y_2,&  \\[1mm]
[Y_1,Y_1]=X_{n-1} & \\ [1mm] [Y_1,Y_2]=X_{n} &
\end{array}\right.
\end{array}
$$
\end{lem}
\begin{rem}All the above Leibniz superalgebras have nilindex
$n+1$ ($=(n+1)+m-2$), two units less than the total dimension of
the Leibniz superalgebras.
\end{rem}
\begin{proof} The family of ${\cal ZF}^{n+1,2}$, with $n\geq 2$, can be
expressed, in an adapted basis $\{X_0,X_1,\dots,X_n, Y_1,Y_2\}$,
by $$\left\{\begin{array}{ll}
[X_i,X_0]=X_{i+1},& 0\leq i \leq n-1 \\[1mm]
[Y_1,X_0]=Y_2 & \\[1mm]
[X_0,Y_1]=a Y_2& \\[1mm]
[Y_1,Y_1]=\displaystyle\sum_{i=0}^{n}b_{11}^{i}X_i & \\[1mm]
[Y_1,Y_2]=\displaystyle\sum_{i=0}^{n}b_{12}^{i}X_i & \\[1mm]
[Y_2,Y_1]=\displaystyle\sum_{i=0}^{n}b_{21}^{i}X_i & \\[1mm]
[Y_2,Y_2]=b_{22}^{n}X_n & \\[1mm]
\end{array}\right.
$$ verifying the graded Leibniz identity. This fact leads to
$$\left\{\begin{array}{ll}
[X_i,X_0]=X_{i+1},& 0\leq i \leq n-1 \\[1mm]
[Y_1,X_0]=Y_2 & \\[1mm]
[X_0,Y_1]=a Y_2& \\[1mm]
[Y_1,Y_1]=b_{11}^{n-1}X_{n-1}+b_{11}^{n}X_{n} & \\[1mm]
[Y_1,Y_2]=b_{12}^{n}X_n & \\[1mm]
[Y_2,Y_1]=b_{21}^{n}X_n & \\[1mm]
\end{array}\right.
$$ with the restrictions: $$\begin{array}{l} ab_{21}^{n}=0\\
(a+1)b_{12}^n=0\\
b_{12}^{n}-b_{11}^{n-1}+b_{21}^{n}=0
\end{array}$$
Now we consider different cases:
\begin{itemize}
\item Case 1. $b_{12}^{n}=0$. In this case we have that
$b_{11}^{n-1}=b_{21}^{n}=\alpha$ and $a\alpha=0$.

\begin{itemize}

\item Case 1.1. $a=0$.

\begin{itemize}
\item Case 1.1.1. $\alpha=0$. As $b_{11}^{n}$ implies that the law
corresponds to a degenerate case, it is necessary that
$b_{11}^{n}\neq 0$. Applying the change of scale $\{X'_0=bX_0,
Y'_1=Y_1, Y'_2=bY_2\}$ taking $b$ as any complex root of the
equation $x^{n+1}-b_{11}^{n}=0$, we have $\mu_1^0$.

\item Case 1.1.2. $\alpha \neq 0$. By using the change of basis
that follows
 $$\left\{\begin{array}{ll} X'_0=bX_0+cX_1&\\
X'_i=b^{i+1}X_i+b^{i}cX_{i+1}, & 1 \leq i \leq n-1\\
X'_n=b^{n+1}X_n & \\
Y'_1=Y_1&\\
Y'_2=bY_2&
\end{array}\right.$$
with $b$ any complex root of the equation $x^{n}-\alpha=0$ and
$c=\displaystyle\frac{b_{11}^n}{b^{n-1}}$, we obtain $\mu_2$.
\end{itemize}

\item Case 1.2. $a \neq 0$. As we have argued before
$b_{11}^{n}\neq 0$ (in other case it will be a degenerate case).
Therefore, without loss of
 generality we may assume that $b_{11}^{n}=1$ (by a simple change of
 scale). If we denote $a$ by $\alpha$, we obtain $\mu_{1}^{\alpha}$ with $\alpha \neq 0$.
\end{itemize}

\item Case 2. $b_{12}^{n}\neq 0$. In this case we have that $a=-1$
and $b_{11}^{n-1}=b_{12}^{n}\neq 0$ and $b_{21}^{n}=0$. By
repeating the change of basis in the case 1.1.2. we obtain
$\mu_{3}$. \end{itemize}

After, by applying a generic change of basis it can be seen that the above
superalgebras are pairwise non-isomorphic.

\end{proof}

\

\begin{thm} The orbits $O(\mu_2),  \ O(\mu_3)$ are Zariski open subsets
of $N^{n+1,2}.$
\end{thm}

\begin{proof} We have that $dimL(\mu_i)=2$ for $i=1,3$, but for $i=2$
this invariant is equal to one. Therefore, $O(\mu_2) \not\subset
\cup_{\alpha\in C} clO(\mu_1^{\alpha})\cup clO(\mu_3)$ which
implies that  $O(\mu_2) \not\subset \cup_{\alpha\in C}
clO(\mu_1^{\alpha})\cup clO(\mu_3)\cup N_{n,2}^{n+1,2}\cup
N_{n+1,1}^{n+1,2}$, that is, $O(\mu_2)$ is a Zariski open subset
of $N^{n+1,2}.$

Since $dimZ(\mu_i)=n+1$ for $i=1,2,$ and $dimZ(\mu_3)=n$, then
$O(\mu_3) \not\subset \cup_{\alpha\in C} clO(\mu_1^{\alpha})\cup
clO(\mu_2)$ which implies that $O(\mu_3) \not\subset
\cup_{\alpha\in C} clO(\mu_1^{\alpha})\cup clO(\mu_2)\cup
N_{n,2}^{n+1,2}\cup N_{n+1,1}^{n+1,2}$, that is, $O(\mu_3)$ is a
Zariski open subset of $N^{n+1,2}.$

\end{proof}

The question now is if there exists any Leibniz superalgebra of
nilindex higher than $n+1$ that is the maximal nilindex for
zero-filiform Leibniz superalgebras. The following theorem is a
consequence of the search of such Leibniz superalgebra in the
filiform Leibniz superalgebras.

\

\begin{thm} Let $L$ be any filiform (non Lie) Leibniz superalgebra
$L=L_0\oplus  L_1$ with $dim(L_0)=n+1$, $n>1$
 and $dim(L_1)=2.$ If
$L$ has nilindex $n+2$, then $n=2$ and it will be isomorphic to
the following Leibniz superalgebra that can be expressed in an
adapted basis $\{X_0,X_1,X_2, Y_1,Y_2\}$ by
$$R^{3,2}=\left\{\begin{array}{l}
[X_1,X_0]=X_2 \\[1mm]
[X_0,X_0]=X_2  \\[1mm]
[X_0,Y_1]=\frac{1}{2} Y_2 \\[1mm]
[X_1,Y_1]=\frac{1}{2} Y_2 \\[1mm]
[Y_1,X_0]=Y_2 \\[1mm]
[Y_1,Y_1]=X_0  \\[1mm]
[Y_2,Y_1]=X_2
\end{array}\right.
$$
\end{thm}

\begin{cor}

\

\qquad $f(3,2)=4$,

\qquad $f(n+1,2)=n+1, $ if $n\geq 3$

\end{cor}

\begin{cor}$${\cal M}^{3,2}\not\subset {\cal ZF}^{3,2}$$
\end{cor}

\begin{cor}$${\cal ZF}^{n+1,2} \subset {\cal M}^{n+1,2}, \quad n\geq 3$$
\end{cor}

\

\noindent {\it Proof of the theorem.} For $n>2$, by using the
graded Leibniz identity it is easy to see that any filiform
Leibniz (no Lie) superalgebra $L=L_0\oplus L_1$ with adapted basis
$\{X_0,X_1,X_2,X_3,\dots,X_n,Y_1,Y_2\}$ verifies that
$[Y_1,Y_1]\in <X_2,\dots,X_n>$, which implies that $L$ will have
always nilindex $n$. For $n=2$, all the cases are as described
below except for one case in which it is possible that $X_0 \in
<[Y_1,Y_1]>$. This exception corresponds to the two-parametric
family
$$\left\{\begin{array}{l}
[X_1,X_0]=X_2 \\[1mm]
[X_0,X_0]=X_2  \\[1mm]
[X_0,Y_1]=\lambda Y_2 \\[1mm]
[X_1,Y_1]=\lambda Y_2 \\[1mm]
[Y_1,X_0]=2\lambda Y_2 \\[1mm]
[Y_1,Y_1]=2\lambda \beta X_0  \\[1mm]
[Y_2,Y_1]=\beta X_2
\end{array}\right.$$
with $\lambda\beta\not = 0$. If we apply the change of scale
$\{X'_0=\displaystyle\frac{1}{2} X_0,
X'_1=\displaystyle\frac{1}{2} X_1,X'_2=\displaystyle
\frac{1}{4}X_2, Y'_1=\displaystyle \frac{1}{2\sqrt{\lambda\beta}}
Y_1,Y'_2=\displaystyle \frac{\sqrt{\lambda}}{2\sqrt{\beta}}Y_2 \}$
to the two-parametric family, we obtain $R$. This concludes the
proof. \hfill{\fbox{}}

\section{Leibniz Superalgebras with three-dimensional odd part}
In this section and in the rest, most of classification proofs are
omitted because of they are very laborious and they do not apport
any new idea.

\subsection{Case $n=1$}
\begin{lem}\label{lema2} Let $L$ be
any Leibniz superalgebra $L \in {\cal ZF}^{2,3}$. Then, it is
isomorphic to one of the following Leibniz superalgebras, pairwise
non-isomorphic, that can be expressed in an adapted basis
$\{X_0,X_1,Y_1,Y_2,Y_3\}$ by
$$\begin{array}{ll}\mu_{1} =\left\{\begin{array}{ll}
[X_0,X_0]=X_{1},& \\[1mm]
[Y_1,X_0]=Y_2 & \\[1mm]
[Y_2,X_0]=Y_3 & \\[1mm]
%[X_0,Y_1]=\lambda Y_2,& \lambda \in \CC\\[1mm]
[Y_1,Y_1]=X_1 &
\end{array}\right.&
\mu_{2} =\left\{\begin{array}{ll}
[X_0,X_0]=X_1& \\[1mm]
[Y_1,X_0]=Y_2 & \\[1mm]
[Y_2,X_0]=Y_3 & \\[1mm]
[X_0,Y_1]=Y_3& \\[1mm]
[Y_1,Y_1]=X_1 & \\[1mm]
\end{array}\right.\\ \\
\mu_{3}=\left\{\begin{array}{ll}
[X_0,X_0]=X_1& \\[1mm]
[Y_1,X_0]=Y_2 & \\[1mm]
[Y_2,X_0]=Y_3 & \\[1mm]
[X_0,Y_1]=-Y_2& \\[1mm]
[X_0,Y_2]=-Y_3& \\[1mm]
[Y_1,Y_1]=X_1 & \\[1mm]
\end{array}\right.&
\mu_{4}=\left\{\begin{array}{ll}
[X_0,X_0]=X_1& \\[1mm]
[Y_1,X_0]=Y_2 & \\[1mm]
[Y_2,X_0]=Y_3 & \\[1mm]
[X_0,Y_1]=-Y_2+Y_3& \\[1mm]
[X_0,Y_2]=-Y_3& \\[1mm]
[Y_1,Y_1]=X_1 & \\[1mm]
\end{array}\right.
\\ \\
%\end{array}$$
%$$\begin{array}{ll}
\mu_{5}=\left\{\begin{array}{ll}
[X_0,X_0]=X_1& \\[1mm]
[Y_1,X_0]=Y_2 & \\[1mm]
[Y_2,X_0]=Y_3 & \\[1mm]
[X_0,Y_1]=-Y_2& \\[1mm]
[X_0,Y_2]=-Y_3& \\[1mm]
[Y_1,Y_3]=-X_1 & \\[1mm]
[Y_2,Y_2]=X_1 & \\[1mm]
[Y_3,Y_1]=-X_1 & \\[1mm]
\end{array}\right.&
\mu_{6}=\left\{\begin{array}{ll}
[X_0,X_0]=X_1& \\[1mm]
[Y_1,X_0]=Y_2 & \\[1mm]
[Y_2,X_0]=Y_3 & \\[1mm]
[X_0,Y_1]=\frac{1}{2}Y_2& \\[1mm]
[X_1,Y_1]=\frac{1}{2}Y_3& \\[1mm]
[Y_1,Y_1]=X_0 & \\[1mm]
[Y_2,Y_1]=X_1 & \\[1mm]
\end{array}\right.\\ \\
\end{array}
$$
\end{lem}

\begin{proof} By using a generic change of basis, along with the graded Leibniz
identity we obtain the lemma. \end{proof}

\begin{rem}
It is easily to see that $\mu_6$ is the superalgebra (in theorem
\ref{ayupov}) for the case $(2,3)$, if we only consider the change
of basis:
$e_1=Y_1,e_2=X_0,e_3=\frac{1}{2}Y_2,e_4=\frac{1}{2}X_1,e_5=Y_3$..
\end{rem}
\begin{prop} $$f(2,3)=5$$
\end{prop}
\begin{proof} Is a corollary of the above lemma.
\end{proof}
\begin{prop}$$ {\cal M}^{2,3}=O(\mu_6)$$
\end{prop}
\begin{proof} Is a corollary of the above lemma.
\end{proof}
\begin{prop} $O(\mu_5), \ O(\mu_6)$ are Zariski open subsets of
$N^{2,3}$.
\end{prop}
\begin{proof} Since $\mu_6$ has maximal index of nilpotency we have that
$O(\mu_6)$ is a Zariski open subset of $N^{2,3}.$

Concerning to $O(\mu_5)$ we have that $dimC^2(\mu_i)=1$ for $1\leq
i\leq 4$  and $dimC^2(\mu_5)=2$, which implies that
$O(\mu_5)\not\subset \cup_{1\leq i\leq 4} clO(\mu_i)$. Since
$dimZ(\mu_6)=3$ and $dimZ(\mu_5)=1$ we have that
$O(\mu_5)\not\subset clO(\mu_6)$. Thus, $O(\mu_5)\not\subset
\cup_{1\leq i\leq 6, i\neq 5} clO(\mu_i),$ that is $O(\mu_5)$ is a
Zariski open subset of $N^{2,3}$.

\end{proof}

\subsection{Case $n=2$.}

\

We are going to give two lemmas which will be useful in this
section and in the remaining ones.

\begin{lem}\label{lema} Let $L=L_0\oplus L_1$ be
any Leibniz superalgebra with s-nilindex $(n+1,3)$ or $(n,3)$.
Then, if we call $\{X_0,X_1,\dots,X_n,Y_1,Y_2,Y_3\}$ an adapted
basis of $L$, it will verify

$$[Y_i,Y_j]\in {\cal C}^{n-6+i+j}(L_0) \qquad \ for \ all \ i,j$$

\end{lem}

\begin{proof} Using the graded Leibniz identity we obtain that

\

\noindent$[Y_3,Y_3]=[[Y_2,Y_3],X_0]=[[Y_3,Y_2],X_0]=[[[Y_1,Y_3],X_0],X_0]=[[[Y_3,Y_1],X_0],X_0]=$

$=1/2[[[Y_2,Y_2],X_0],X_0]=1/3[[[[Y_1,Y_2],X_0],X_0],X_0]=1/3[[[[Y_2,Y_1],X_0],X_0],X_0]=$

$=1/6[[[[[Y_1,Y_1],X_0],X_0],X_0],X_0]\in {\cal C}^n(L_0)=
Z_{L_0}(L_0)$

\

which leads to the lemma. \end{proof}

\begin{lem}\label{lema1} Let $L=L_0\oplus L_1$ be
any Leibniz superalgebra with s-nilindex $(n+1,3)$ and $n \geq 2$.
Then, if we denote $\{X_0,X_1,\dots,X_n,Y_1,Y_2,Y_3\}$ an adapted
basis of $L$, it will verify

$$[Y_i,Y_j]\in {\cal C}^{k-2}(L_0) \qquad i+j=k, \quad 3 \leq k \leq 5$$
$$[Y_3,Y_3]\in {\cal C}^{4}(L_0) \cap Z(L_0)$$
\end{lem}
\begin{proof} Using the graded Leibniz identity we obtain that
$$[Y_1,Y_2]+[Y_2,Y_1]=[[Y_1,Y_1],X_0] \in {\cal C}^1(L_0)$$
Furthermore $[Y_1,Y_2], \ [Y_2,Y_1] \in {\cal C}^1(L_0)$. In fact,
if $b_{12}^0 \not = 0$ then by the above equality we have that
$b_{12}^0 = - b_{21}^0$ ($b_{12}^0$ and $b_{21}^0$ are
respectively the coefficients of $X_0$ in $[Y_1,Y_2]$ and
$[Y_2,Y_1]$). By considering the products
$$[X_1,[Y_1,Y_2]]=[[X_1,Y_1],Y_2]+[[X_1,Y_2],Y_1]=b_{12}^0X_2+\gamma
X_3+\dots $$
$$[X_1,[Y_2,Y_1]]=[[X_1,Y_2],Y_1]+[[X_1,Y_1],Y_2]=-b_{12}^0X_2+\beta
X_3+\dots $$ we obtain that $b_{12}^0=0$ which is a contradiction.
Thus, $[Y_1,Y_2], \ [Y_2,Y_1] \in {\cal
C}^1(L_0)$.

Analogously, it can be proved that $[Y_2,Y_2],[Y_1,Y_3],[Y_3,Y_1]
\in {\cal C}^2(L_0)$. Finally by considering
$[Y_3,Y_2]=[[Y_3,Y_1],X_0]$, $[Y_2,Y_3]=[[Y_1,Y_3],X_0]$ and
$[Y_3,Y_3]=[[Y_3,Y_2],X_0]$ we have the lemma. \end{proof}

\

\begin{lem}Let $L$ be
any Leibniz superalgebra $L \in {\cal ZF}^{3,3}$. Then it is
isomorphic to one of the following Leibniz superalgebras, pairwise
non-isomorphic, that can be expressed in an adapted basis
$\{X_0,X_1,X_2,Y_1,Y_2,Y_3\}$ by
$$\begin{array}{ll}
%\mu_1 =\left\{\begin{array}{ll}
%[X_i,X_0]=X_{i+1},& 0\leq i \leq 1\\[1mm]
%[Y_j,X_0]=Y_{j+1} & 1\leq j \leq 2\\[1mm]
%[Y_1,Y_1]=X_2 & \\[1mm]
%\end{array}\right.&
\mu_1 =\left\{\begin{array}{ll}
[X_i,X_0]=X_{i+1},& 0\leq i \leq 1\\[1mm]
[Y_j,X_0]=Y_{j+1} & 1\leq j \leq 2\\[1mm]
%[Y_1,Y_1]=\alpha X_2 & \alpha \in \{0,1\}\\[1mm]
[Y_1,Y_2]=X_2 & \\[1mm]
[Y_2,Y_1]=-X_2 & \\[1mm]
\end{array}\right.&
\mu_2 =\left\{\begin{array}{ll}
[X_i,X_0]=X_{i+1},& 0\leq i \leq 1\\[1mm]
[Y_j,X_0]=Y_{j+1} & 1\leq j \leq 2\\[1mm]
[Y_1,Y_1]=X_2 & \\[1mm]
[Y_1,Y_2]=X_2 & \\[1mm]
[Y_2,Y_1]=-X_2 & \\[1mm]
\end{array}\right.
%\end{array}$$
%$$\begin{array}{ll}
\\[22mm]
\mu_3^{\alpha} =\left\{\begin{array}{ll}
[X_i,X_0]=X_{i+1},& 0\leq i \leq 1\\[1mm]
[Y_j,X_0]=Y_{j+1} & 1\leq j \leq 2\\[1mm]
[Y_1,Y_1]=X_1 & \\[1mm]
[Y_1,Y_2]=\alpha X_2 & \alpha \in \CC\\[1mm]
[Y_2,Y_1]=(1-\alpha)X_2 & \\[1mm]
\end{array}\right.&
\mu_4 =\left\{\begin{array}{ll}
[X_i,X_0]=X_{i+1},& 0\leq i \leq 1\\[1mm]
[Y_j,X_0]=Y_{j+1} & 1\leq j \leq 2\\[1mm]
[X_0,Y_1]=Y_3& \\[1mm]
[Y_1,Y_1]=X_2 & \\[1mm]
\end{array}\right. \\ [22mm]
\mu_5^{\alpha} =\left\{\begin{array}{ll}
[X_i,X_0]=X_{i+1},& 0\leq i \leq 1\\[1mm]
[Y_j,X_0]=Y_{j+1} & 1\leq j \leq 2\\[1mm]
[X_0,Y_1]=Y_3& \\[1mm]
[Y_1,Y_1]=\alpha X_2 & \alpha \in \CC\\[1mm]
[Y_1,Y_2]=X_2 & \\[1mm]
[Y_2,Y_1]=-X_2 & \\[1mm]
\end{array}\right. &
\mu_6^{\alpha} =\left\{\begin{array}{ll}
[X_i,X_0]=X_{i+1},& 0\leq i \leq 1\\[1mm]
[Y_j,X_0]=Y_{j+1} & 1\leq j \leq 2\\[1mm]
[X_0,Y_1]=Y_3& \\[1mm]
[Y_1,Y_1]=X_1 & \\[1mm]
[Y_1,Y_2]=\alpha X_2 & \alpha \in \CC\\[1mm]
[Y_2,Y_1]=(1-\alpha)X_2 & \\[1mm]
\end{array}\right. \\ [22mm]
\mu_{7}^{\a} =\left\{\begin{array}{ll}
[X_i,X_0]=X_{i+1},& 0\leq i \leq 1\\[1mm]
[Y_j,X_0]=Y_{j+1} & 1\leq j \leq 2\\[1mm]
[X_0,Y_1]=\a Y_2& \a \in \CC\\[1mm]
[X_1,Y_1]=\a Y_3& \\[1mm]
[Y_1,Y_1]=X_2 & \\[1mm]
\end{array}\right. &
\mu_{8}=\left\{\begin{array}{ll}
[X_i,X_0]=X_{i+1},& 0\leq i \leq 1\\[1mm]
[Y_j,X_0]=Y_{j+1} & 1\leq j \leq 2\\[1mm]
%[X_0,Y_1]=-Y_2+\a Y_3& \a \in \{0,1\} \\[1mm]
[X_0,Y_2]=-Y_3& \\[1mm]
[Y_1,Y_1]=X_1 & \\[1mm]
[Y_1,Y_2]=X_2 & \\[1mm]
\end{array}\right. \\ [22mm]
\mu_{9}=\left\{\begin{array}{ll}
[X_i,X_0]=X_{i+1},& 0\leq i \leq 1\\[1mm]
[Y_j,X_0]=Y_{j+1} & 1\leq j \leq 2\\[1mm]
[X_0,Y_1]=-Y_2+Y_3&  \\[1mm]
[X_0,Y_2]=-Y_3& \\[1mm]
[Y_1,Y_1]=X_1 & \\[1mm]
[Y_1,Y_2]=X_2 & \\[1mm]
\end{array}\right. &
\mu_{10}=\left\{\begin{array}{ll}
[X_i,X_0]=X_{i+1},& 0\leq i \leq 1\\[1mm]
[Y_j,X_0]=Y_{j+1} & 1\leq j \leq 2\\[1mm]
[X_0,Y_1]=-Y_2& \\[1mm]
[X_0,Y_2]=-Y_3& \\[1mm]
[Y_1,Y_3]=-X_2 & \\[1mm]
[Y_2,Y_2]=X_2 & \\[1mm]
[Y_3,Y_1]=-X_2 & \\[1mm]
\end{array}\right.
% \\ \\
\end{array}$$
$$\begin{array}{ll}
\mu_{11}=\left\{\begin{array}{ll}
[X_i,X_0]=X_{i+1},& 0\leq i \leq 1\\[1mm]
[Y_j,X_0]=Y_{j+1} & 1\leq j \leq 2\\[1mm]
[X_0,Y_1]=-Y_2& \\[1mm]
[X_0,Y_2]=-Y_3& \\[1mm]
[Y_1,Y_1]=X_1 & \\[1mm]
[Y_1,Y_2]=X_2 & \\[1mm]
[Y_1,Y_3]=-X_2 & \\[1mm]
[Y_2,Y_2]=X_2 & \\[1mm]
[Y_3,Y_1]=-X_2 & \\[1mm]
\end{array}\right. &
 \mu_{12}=\left\{\begin{array}{ll}
[X_i,X_0]=X_{i+1},& 0\leq i \leq 1\\[1mm]
[Y_j,X_0]=Y_{j+1} & 1\leq j \leq 2\\[1mm]
[X_0,Y_1]=\frac{1}{2}Y_2& \\[1mm]
[X_1,Y_1]=\frac{1}{2}Y_3& \\[1mm]
[Y_1,Y_1]=X_0 & \\[1mm]
[Y_2,Y_1]=X_1 & \\[1mm]
[Y_3,Y_1]=X_2 & \\[1mm]
\end{array}\right. \\ \\
\end{array}
$$
\end{lem}
\begin{proof} By using a generic change of basis, along with the graded Leibniz
identity and the lemma \ref{lema1}, we obtain the result.
\end{proof}

\begin{rem}
It is easily to see that $\mu_{12}$ is the Leibniz superalgebra of
the theorem \ref{ayupov} for the case $(3,3)$.
\end{rem}
\begin{prop} $$f(3,3)=6$$
\end{prop}
\begin{proof} Is a corollary of the above lemma.
\end{proof}
\begin{prop}$$ {\cal M}^{3,3}=O(\mu_{12})$$
\end{prop}
\begin{proof} Is a corollary of the above lemma.
\end{proof}
\begin{prop} $O(\mu_{11}),$ $O(\mu_{12})$ are Zariski open subsets of
$N^{3,3}.$
\end{prop}

\begin{proof} Since $dimC^3(\mu_i)=0$ for $1\leq i\leq 9$ and
$dimC^3(\mu_{10})=dimC^3(\mu_{11})=1$ and $dimC^3(\mu_{12})=3$, we
have that $O(\mu_{12})$ is an open subset of $N^{3,3}$ and
$O(\mu_{11})\not\subset \cup_{1\leq i\leq 9}clO(\mu_i).$ Note that
$O(\mu_{10})\subset clO(\mu_{11}),$ in fact if we take
$f_t=f_{to}+f_{t1},$ where $f_{to}(X_0)=t^{-1}X_0, \
f_{to}(X_1)=t^{-2}X_1, \ f_{to}(X_2)=t^{-3}X_2,$ and
$f_{t1}(Y_1)=t^{-\frac{1}{2}}Y_1, \
f_{t1}(Y_2)=t^{-\frac{3}{2}}Y_2, \
f_{t1}(Y_3)=t^{-\frac{5}{2}}Y_3,$ we obtain that if $t\rightarrow
0$ then $O(\mu_{10})\subset clO(\mu_{11}).$

Since $dimCent(\mu_{12})=4$ and $dimCent(\mu_{11})=2$, we have
that $O(\mu_{11})\not\subset clO(\mu_{12}).$ Thus
$O(\mu_{11})\not\subset \cup_{{1\leq i\leq 12}, i\neq
11}clO(\mu_{i})$, that is, $O(\mu_{11})$ is an open subset of
$N^{3,3}$
\end{proof}

\subsection{Case $n=3$.}

\

\begin{lem}Let $L$ be
any Leibniz superalgebra $L \in {\cal ZF}^{4,3}$. Then, it is
isomorphic to one of the following Leibniz superalgebras that can
be expressed in an adapted basis $\{X_0,X_1,X_2,X_3,Y_1,Y_2,Y_3\}$
by

$$\begin{array}{ll}\mu_1 =\left\{\begin{array}{ll}
[X_i,X_0]=X_{i+1},& 0\leq i \leq 2\\[1mm]
[Y_j,X_0]=Y_{j+1} & 1\leq j \leq 2\\[1mm]
[X_0,Y_1]=-Y_2 +Y_3& \\[1mm]
[X_0,Y_2]=-Y_3& \\[1mm]
[Y_1,Y_1]=X_3 & \\[1mm]
\end{array}\right.&
\mu_2^{\a} =\left\{\begin{array}{ll}
[X_i,X_0]=X_{i+1},& 0\leq i \leq 2\\[1mm]
[Y_j,X_0]=Y_{j+1} & 1\leq j \leq 2\\[1mm]
[X_0,Y_1]=-Y_2 & \\[1mm]
[X_0,Y_2]=(-1-\a)Y_3& \a \in \CC\\[1mm]
[X_1,Y_1]=\a Y_3& \\[1mm]
[Y_1,Y_1]=X_3 & \\[1mm]
\end{array}\right.\\ [25mm]
\mu_3^{\a} =\left\{\begin{array}{ll}
[X_i,X_0]=X_{i+1},& 0\leq i \leq 2\\[1mm]
[Y_j,X_0]=Y_{j+1} & 1\leq j \leq 2\\[1mm]
[X_0,Y_1]=-Y_2 & \\[1mm]
[X_0,Y_2]=(-1-\a)Y_3& \a \in \CC\\[1mm]
[X_1,Y_1]=\a Y_3& \\[1mm]
[Y_1,Y_1]=X_2& \\[1mm]
[Y_1,Y_2]=X_3 & \\[1mm]
\end{array}\right.&
\mu_4^{\a} =\left\{\begin{array}{ll}
[X_i,X_0]=X_{i+1},& 0\leq i \leq 2\\[1mm]
[Y_j,X_0]=Y_{j+1} & 1\leq j \leq 2\\[1mm]
[X_0,Y_1]=-Y_2 +Y_3& \\[1mm]
[X_0,Y_2]=(-1-\a)Y_3& \a \in \CC\\[1mm]
[X_1,Y_1]=\a Y_3& \\[1mm]
[Y_1,Y_1]=X_2 & \\[1mm]
[Y_1,Y_2]=X_3 & \\[1mm]
\end{array}\right.
%\end{array}$$
%$$\begin{array}{ll}
\\ [25mm]
\mu_5=\left\{\begin{array}{ll}
[X_i,X_0]=X_{i+1},& 0\leq i \leq 2\\[1mm]
[Y_j,X_0]=Y_{j+1} & 1\leq j \leq 2\\[1mm]
[X_0,Y_1]=Y_3& \\[1mm]
[Y_1,Y_1]=X_3 &\\[1mm]
\end{array}\right.&
\mu_6=\left\{\begin{array}{ll}
[X_i,X_0]=X_{i+1},& 0\leq i \leq 2\\[1mm]
[Y_j,X_0]=Y_{j+1} & 1\leq j \leq 2\\[1mm]
%[X_0,Y_1]=\a Y_3& \a \in \{0,1\} \\[1mm]
[Y_1,Y_1]=X_2 &\\[1mm]
[Y_2,Y_1]=X_3 &\\[1mm]
\end{array}\right.\\ [25mm]
\mu_7=\left\{\begin{array}{ll}
[X_i,X_0]=X_{i+1},& 0\leq i \leq 2\\[1mm]
[Y_j,X_0]=Y_{j+1} & 1\leq j \leq 2\\[1mm]
[X_0,Y_1]= Y_3&  \\[1mm]
[Y_1,Y_1]=X_2 &\\[1mm]
[Y_2,Y_1]=X_3 &\\[1mm]
\end{array}\right.&
\mu_8^{\a} =\left\{\begin{array}{ll}
[X_i,X_0]=X_{i+1},& 0\leq i \leq 2\\[1mm]
[Y_j,X_0]=Y_{j+1} & 1\leq j \leq 2\\[1mm]
[X_0,Y_1]=\a Y_2  & \a \in \CC\\[1mm]
[X_1,Y_1]=\a Y_3& \\[1mm]
[Y_1,Y_1]= X_3 & \\[1mm]
\end{array}\right.\\ [25mm]
\mu_9=\left\{\begin{array}{ll}
[X_i,X_0]=X_{i+1},& 0\leq i \leq 2\\[1mm]
[Y_j,X_0]=Y_{j+1} & 1\leq j \leq 2\\[1mm]
[X_0,Y_1]=-Y_2 & \\[1mm]
[X_0,Y_2]=-Y_3& \\[1mm]
%[Y_1,Y_1]=\a X_2 & \a \in \{0,1\} \\[1mm]
%[Y_1,Y_2]=\a X_3 & \\[1mm]
[Y_1,Y_3]=X_3&\\[1mm]
[Y_2,Y_2]=-X_3&\\[1mm]
[Y_3,Y_1]=X_3&
\end{array}\right.&
\mu_{10}=\left\{\begin{array}{ll}
[X_i,X_0]=X_{i+1},& 0\leq i \leq 2\\[1mm]
[Y_j,X_0]=Y_{j+1} & 1\leq j \leq 2\\[1mm]
[X_0,Y_1]=-Y_2 & \\[1mm]
[X_0,Y_2]=-Y_3& \\[1mm]
[Y_1,Y_1]=X_2 &  \\[1mm]
[Y_1,Y_2]= X_3 & \\[1mm]
[Y_1,Y_3]=X_3&\\[1mm]
[Y_2,Y_2]=-X_3&\\[1mm]
[Y_3,Y_1]=X_3&
\end{array}\right.
\end{array}$$
$$\begin{array}{ll}
\mu_{11}^{\a}=\left\{\begin{array}{ll}
[X_i,X_0]=X_{i+1},& 0\leq i \leq 2\\[1mm]
[Y_j,X_0]=Y_{j+1} & 1\leq j \leq 2\\[1mm]
[X_0,Y_1]=-Y_2+ Y_3 & \\[1mm]
[X_0,Y_2]=-Y_3& \\[1mm]
[Y_1,Y_1]=\a X_2 & \a \in \CC\\[1mm]
[Y_1,Y_2]=\a X_3 & \\[1mm]
[Y_1,Y_3]=X_3&\\[1mm]
[Y_2,Y_2]=-X_3&\\[1mm]
[Y_3,Y_1]=X_3&
\end{array}\right.&
\mu_{12}=\left\{\begin{array}{ll}
[X_i,X_0]=X_{i+1},& 0\leq i \leq 2\\[1mm]
[Y_j,X_0]=Y_{j+1} & 1\leq j \leq 2\\[1mm]
%[X_0,Y_1]=\a Y_3 & \a \in \{0,1\} \\[1mm]
[Y_1,Y_1]=X_1 & \\[1mm]
[Y_2,Y_1]=X_2 & \\[1mm]
[Y_3,Y_1]=X_3&\\
\end{array}\right.\\ \\
\mu_{13}=\left\{\begin{array}{ll}
[X_i,X_0]=X_{i+1},& 0\leq i \leq 2\\[1mm]
[Y_j,X_0]=Y_{j+1} & 1\leq j \leq 2\\[1mm]
[X_0,Y_1]=Y_3 & \\[1mm]
[Y_1,Y_1]=X_1 & \\[1mm]
[Y_2,Y_1]=X_2 & \\[1mm]
[Y_3,Y_1]=X_3&\\
\end{array}\right. &
\end{array}$$
\end{lem}

\begin{proof} By using a generic change of basis, along with the graded Leibniz
identity we obtain the result. \end{proof}

\begin{rem} All the above zero-filiform Leibniz algebras have
nilindex $4$, three units less than the total dimension of the
superalgebra.
\end{rem}

\begin{prop} There exists one non empty subset J of $\CC$ such that the set $cl(\cup_{\alpha\in
J}O(\mu_{11}^{\alpha}))$ is an irreducible component of $N^{4,3}$
\end{prop}
\begin{proof} For $1\leq i\leq 13, i\neq {9,10,11}$ the invariant
$dimZ(\mu_i)$ is higher or equal to 4, but for $i={9,10,11}$ it is
less than 4. And $\mu_9\subset clO(\mu_{10}), \mu_{10}\subset
clO(\mu_{11}^1),$ in fact, if we take the following maps

$f_t=f_{to}+f_{t1},$ where $f_{to}(X_0)=tX_0,$ \
$f_{to}(X_1)=t^2X_1,$ \ $f_{to}(X_2)=t^3X_2,$ \
$f_{to}(X_3)=t^4X_3,$ and $f_{t1}(Y_1)=tY_1,$ \
$f_{t1}(Y_2)=t^2Y_2,$ \ $f_{t1}(Y_3)=t^3Y_3,$

and correspondingly

$g_t=g_{to}+g_{t1},$ where $g_{to}(X_0)=tX_0,$ \
$g_{to}(X_1)=t^2X_1,$ \ $g_{to}(X_2)=t^3X_2,$ \
$g_{to}(X_3)=t^4X_3,$ and $g_{t1}(Y_1)=tY_1,$ \
$g_{t1}(Y_2)=t^2Y_2,$ \ $g_{t1}(Y_3)=t^3Y_3,$

we obtain that if $t\rightarrow 0$  then $\mu_9\subset
clO(\mu_{10}), \ \mu_{10}\subset clO(\mu_{11}^1).$

We have that exists such subset $J$ of $\CC$ stated in the
proposition.
\end{proof}

The question now is if there exists any Leibniz superalgebra of
nilindex higher than $4$ that is the maximal nilindex maximal for
the zero-filiform Leibniz superalgebras. When we search this
superalgebra in filiform Leibniz superalgebras (the only remain
family of Leibniz superalgebras that can arise nilindex $6$, one
unit smaller than the total dimension) appears the following
proposition.

\begin{prop} \label{prop4,3} Let $L$ be
any filiform (non Lie) Leibniz superalgebra $L=L_0\oplus L_1$ with
$dim(L_0)=4$,
 and $dim(L_1)=3.$ If
$L$ has nilindex $6$, then it will be isomorphic to the following
Leibniz superalgebra that can be expressed in an adapted basis
$\{X_0,X_1,X_2,X_3, Y_1,Y_2,Y_3\}$ by
$$R^{4,3}=\left\{\begin{array}{l}
[X_1,X_0]=X_2 \\[1mm]
[X_2,X_0]=X_3\\[1mm]
[X_0,X_0]=X_2  \\[1mm]
[X_0,Y_1]=\frac{1}{2} Y_2 \\[1mm]
[X_1,Y_1]=\frac{1}{2} Y_2 \\[1mm]
[X_2,Y_1]=\frac{1}{2} Y_3 \\[1mm]
[Y_1,X_0]=Y_2 \\[1mm]
[Y_2,X_0]=Y_3 \\[1mm]
[Y_1,Y_1]=X_0  \\[1mm]
[Y_2,Y_1]=X_2 \\[1mm]
[Y_3,Y_1]=X_3\\
\end{array}\right.
$$
\end{prop}
\begin{proof} By using lemma \ref{lema} and  graded Leibniz identity we
obtain that the only family of laws with nilindex $6$ is
$$\left\{\begin{array}{l}
[X_1,X_0]=X_2 \\[1mm]
[X_2,X_0]=X_3\\[1mm]
[X_0,X_0]=X_2  \\[1mm]
[X_0,Y_1]=\displaystyle \frac{b_{11}^0}{2(b_{11}^0+b_{11}^1)} Y_2
-
\displaystyle \frac{b_{11}^2 b_{11}^0}{2(b_{11}^0+b_{11}^1)^2}Y_3\\[1mm]
[X_1,Y_1]= \displaystyle \frac{b_{11}^0}{2(b_{11}^0+b_{11}^1)} Y_2
-\displaystyle \frac{b_{11}^2 b_{11}^0}{2(b_{11}^0+b_{11}^1)^2}Y_3\\[1mm]
[X_2,Y_1]= \displaystyle \frac{b_{11}^0}{2(b_{11}^0+b_{11}^1)}Y_3 \\[1mm]
[Y_1,X_0]=Y_2 \\[1mm]
[Y_2,X_0]=Y_3 \\[1mm]
[Y_1,Y_1]=b_{11}^0X_0+b_{11}^1 X_1+b_{11}^2X_2+b_{11}^3 X_3 \\[1mm]
[Y_2,Y_1]=(b_{11}^0+b_{11}^1)X_2 + b_{11}^2X_3\\[1mm]
[Y_3,Y_1]=(b_{11}^0+b_{11}^1)X_3\\
\end{array}\right.
$$
with the restrictions $b_{11}^0 \not = 0$ and $b_{11}^1 \not = -
b_{11}^0$. By applying the change of basis
$$\left\{\begin{array}{l}
X'_0=b_{11}^0X_0+b_{11}^1 X_1+b_{11}^2X_2+b_{11}^3 X_3 \\[1mm]
X'_1=(b_{11}^0+b_{11}^1) X_1+b_{11}^2X_2+b_{11}^3 X_3\\[1mm]
X'_2=b_{11}^0(b_{11}^0+b_{11}^1) X_2+b_{11}^0 b_{11}^2X_3\\[1mm]
X'_3=(b_{11}^0)^2(b_{11}^0+b_{11}^1) X_3\\[1mm]
Y'_1=Y_1\\[1mm]
Y'_2=b_{11}^0Y_2\\[1mm]
Y'_3=(b_{11}^0)^2 Y_3
\end{array}\right.
$$
$R^{4,3}$ is obtained. \end{proof}

\begin{prop} $$f(4,3)=6$$
\end{prop}
\begin{proof} Is a corollary of the above proposition.
\end{proof}
\begin{prop}$$ {\cal M}^{4,3}\not \subset {\cal ZF}^{4,3}$$
\end{prop}
\begin{proof} Is a corollary of the proposition \ref{prop4,3}.
\end{proof}
%\begin{cor} ${\cal O}(R^{4,3})$ is a Zariski open subset of
%${\cal N}^{4,3}$.
%\end{cor}

\subsection{Case $n\geq 4$.}
\begin{lem}Let $L$ be
any Leibniz superalgebra $L \in {\cal ZF}^{n+1,3}$, with $n\geq
4$. Then it is isomorphic to one of the following Leibniz
superalgebras, pairwise non-isomorphic, that can be expressed in
an adapted basis $\{X_0,X_1,X_2,X_3,X_4,\dots,X_n,Y_1,Y_2,Y_3\}$
by $$\begin{array}{ll}\mu_1 =\left\{\begin{array}{ll}
[X_i,X_0]=X_{i+1},& 0\leq i \leq n-1\\[1mm]
[Y_j,X_0]=Y_{j+1} & 1\leq j \leq 2\\[1mm]
[X_0,Y_1]=-Y_2 +Y_3& \\[1mm]
[X_0,Y_2]=-Y_3& \\[1mm]
[Y_1,Y_1]=X_n & \\[1mm]
\end{array}\right.&
\mu_2^{\a} =\left\{\begin{array}{ll}
[X_i,X_0]=X_{i+1},& 0\leq i \leq n-1\\[1mm]
[Y_j,X_0]=Y_{j+1} & 1\leq j \leq 2\\[1mm]
[X_0,Y_1]=-Y_2 & \\[1mm]
[X_0,Y_2]=(-1-\a)Y_3& \a \in \CC\\[1mm]
[X_1,Y_1]=\a Y_3& \\[1mm]
[Y_1,Y_1]=X_n & \\[1mm]
\end{array}\right.\\ [20mm]
%\end{array}$$
%$$\begin{array}{ll}
\mu_3^{\a} =\left\{\begin{array}{ll}
[X_i,X_0]=X_{i+1},& 0\leq i \leq n-1\\[1mm]
[Y_j,X_0]=Y_{j+1} & 1\leq j \leq 2\\[1mm]
[X_0,Y_1]=-Y_2 & \\[1mm]
[X_0,Y_2]=(-1-\a)Y_3& \a \in \CC\\[1mm]
[X_1,Y_1]=\a Y_3& \\[1mm]
[Y_1,Y_1]=X_{n-1}& \\[1mm]
[Y_1,Y_2]=X_n & \\[1mm]
\end{array}\right.&
\mu_4^{\a} =\left\{\begin{array}{ll}
[X_i,X_0]=X_{i+1},& 0\leq i \leq n-1\\[1mm]
[Y_j,X_0]=Y_{j+1} & 1\leq j \leq 2\\[1mm]
[X_0,Y_1]=-Y_2 +Y_3& \\[1mm]
[X_0,Y_2]=(-1-\a)Y_3& \a \in \CC\\[1mm]
[X_1,Y_1]=\a Y_3& \\[1mm]
[Y_1,Y_1]=X_{n-1} & \\[1mm]
[Y_1,Y_2]=X_n & \\[1mm]
\end{array}\right.\\ [20mm]
%\end{array}$$
%$$\begin{array}{ll}
 \mu_5=\left\{\begin{array}{ll}
[X_i,X_0]=X_{i+1},& 0\leq i \leq n-1\\[1mm]
[Y_j,X_0]=Y_{j+1} & 1\leq j \leq 2\\[1mm]
[X_0,Y_1]=Y_3&  \\[1mm]
[Y_1,Y_1]=X_n &\\[1mm]
\end{array}\right.&
\mu_6=\left\{\begin{array}{ll}
[X_i,X_0]=X_{i+1},& 0\leq i \leq n-1\\[1mm]
[Y_j,X_0]=Y_{j+1} & 1\leq j \leq 2\\[1mm]
%[X_0,Y_1]=\a Y_3& \a \in \{0,1\} \\[1mm]
[Y_1,Y_1]=X_{n-1} &\\[1mm]
[Y_2,Y_1]=X_n &\\[1mm]
\end{array}\right.\\ [20mm]
%\end{array}$$
%$$\begin{array}{ll}
\mu_7=\left\{\begin{array}{ll}
[X_i,X_0]=X_{i+1},& 0\leq i \leq n-1\\[1mm]
[Y_j,X_0]=Y_{j+1} & 1\leq j \leq 2\\[1mm]
[X_0,Y_1]=Y_3&  \\[1mm]
[Y_1,Y_1]=X_{n-1} &\\[1mm]
[Y_2,Y_1]=X_n &\\[1mm]
\end{array}\right.&
\mu_8^{\a} =\left\{\begin{array}{ll}
[X_i,X_0]=X_{i+1},& 0\leq i \leq n-1\\[1mm]
[Y_j,X_0]=Y_{j+1} & 1\leq j \leq 2\\[1mm]
[X_0,Y_1]=\a Y_2  & \a \in \CC\\[1mm]
[X_1,Y_1]=\a Y_3& \\[1mm]
[Y_1,Y_1]= X_n & \\[1mm]
\end{array}\right.
\end{array}$$
$$\begin{array}{ll}
\mu_9=\left\{\begin{array}{ll}
[X_i,X_0]=X_{i+1},& 0\leq i \leq n-1\\[1mm]
[Y_j,X_0]=Y_{j+1} & 1\leq j \leq 2\\[1mm]
[X_0,Y_1]=-Y_2 & \\[1mm]
[X_0,Y_2]=-Y_3& \\[1mm]
%[Y_1,Y_1]=\a X_{n-1} & \a \in \{0,1\} \\[1mm]
%[Y_1,Y_2]=\a X_n & \\[1mm]
[Y_1,Y_3]=X_n&\\[1mm]
[Y_2,Y_2]=-X_n&\\[1mm]
[Y_3,Y_1]=X_n&
\end{array}\right.&
\mu_{10}=\left\{\begin{array}{ll}
[X_i,X_0]=X_{i+1},& 0\leq i \leq n-1\\[1mm]
[Y_j,X_0]=Y_{j+1} & 1\leq j \leq 2\\[1mm]
[X_0,Y_1]=-Y_2 & \\[1mm]
[X_0,Y_2]=-Y_3& \\[1mm]
[Y_1,Y_1]=X_{n-1} &  \\[1mm]
[Y_1,Y_2]= X_n & \\[1mm]
[Y_1,Y_3]=X_n&\\[1mm]
[Y_2,Y_2]=-X_n&\\[1mm]
[Y_3,Y_1]=X_n&
\end{array}\right.
\end{array}$$

$$\begin{array}{ll}
\mu_{11}^{\a}=\left\{\begin{array}{ll}
[X_i,X_0]=X_{i+1},& 0\leq i \leq n-1\\[1mm]
[Y_j,X_0]=Y_{j+1} & 1\leq j \leq 2\\[1mm]
[X_0,Y_1]=-Y_2+ Y_3 & \\[1mm]
[X_0,Y_2]=-Y_3& \\[1mm]
[Y_1,Y_1]=\a X_{n-1} & \a \in \CC\\[1mm]
[Y_1,Y_2]=\a X_n & \\[1mm]
[Y_1,Y_3]=X_n&\\[1mm]
[Y_2,Y_2]=-X_n&\\[1mm]
[Y_3,Y_1]=X_n&
\end{array}\right.&
\mu_{12}=\left\{\begin{array}{ll}
[X_i,X_0]=X_{i+1},& 0\leq i \leq n-1\\[1mm]
[Y_j,X_0]=Y_{j+1} & 1\leq j \leq 2\\[1mm]
%[X_0,Y_1]=\a Y_3 & \a \in \{0,1\} \\[1mm]
[Y_1,Y_1]=X_{n-2} & \\[1mm]
[Y_2,Y_1]=X_{n-1} & \\[1mm]
[Y_3,Y_1]=X_n&\\
\end{array}\right.\\
\mu_{13}=\left\{\begin{array}{ll}
[X_i,X_0]=X_{i+1},& 0\leq i \leq n-1\\[1mm]
[Y_j,X_0]=Y_{j+1} & 1\leq j \leq 2\\[1mm]
[X_0,Y_1]=Y_3 &  \\[1mm]
[Y_1,Y_1]=X_{n-2} & \\[1mm]
[Y_2,Y_1]=X_{n-1} & \\[1mm]
[Y_3,Y_1]=X_n&\\
\end{array}\right.&
\mu_{14}^{\a}=\left\{\begin{array}{ll}
[X_i,X_0]=X_{i+1},& 0\leq i \leq n-1\\[1mm]
[Y_j,X_0]=Y_{j+1} & 1\leq j \leq 2\\[1mm]
[X_0,Y_1]=-Y_2+\a Y_3 & \a \in \CC \\[1mm]
[X_0,Y_2]=-Y_3 & \\[1mm]
[Y_1,Y_1]=X_{n-2} & \\[1mm]
[Y_1,Y_2]=X_{n-1}+\a X_n & \\[1mm]
[Y_2,Y_1]=\a X_{n} & \\[1mm]
[Y_2,Y_2]=X_{n} & \\[1mm]
[Y_3,Y_1]=-X_{n}&\\
\end{array}\right.\\ \\
\mu_{15}^{\a}=\left\{\begin{array}{ll}
[X_i,X_0]=X_{i+1},& 0\leq i \leq n-1\\[1mm]
[Y_j,X_0]=Y_{j+1} & 1\leq j \leq 2\\[1mm]
[X_0,Y_1]=- Y_2 & \\[1mm]
[X_0,Y_2]=- Y_3 & \\[1mm]
[Y_1,Y_1]=X_{n-2} & \\[1mm]
[Y_1,Y_2]=X_{n-1} & \\[1mm]
[Y_1,Y_3]=\a X_n & \a \in \CC-\{0\}\\ [1mm]
[Y_2,Y_2]=(1-\a )X_{n}&\\[1mm]
[Y_3,Y_1]=(\a -1)X_{n}&\\[1mm]
\end{array}\right.
\end{array}$$
\end{lem}

\begin{proof} By using a generic change of basis, along with the graded Leibniz
identity we obtain the result. \end{proof}

\begin{rem} All the above zero-filiform Leibniz algebras have
nilindex $n+1$, three units smaller than the total dimension of the
superalgebra.
\end{rem}

\

 \begin{thm}  Let $L$ be an arbitrary (non Lie) Leibniz
filiform superalgebra with $gz(L)=(n,1|3)$ and $n\geq 4$. Then $L$
always has nilindex $n$.
\end{thm}

\begin{proof} Is a consequence of the above lemma, of the lemma \ref{lema}, and of
graded Leibniz identity.
\end{proof}

\

\begin{thm}  $$ {\cal M}^{n+1,3} \not \subset {\cal
ZF}^{n+1,3},\qquad \mbox{ if } n=4$$
\end{thm}

\begin{proof} Consider $R^{4,3}\oplus \{X_4\} \notin {\cal ZF}^{5,3}$ with
nilindex $5$. \end{proof}

\noindent{\bf Conjecture 1.}
$$f(n+1,3)=n+1, \qquad  \mbox{ if } n\geq 4$$
$$ {\cal M}^{n+1,3} \subset {\cal
ZF}^{n+1,3},\qquad \mbox{ if } n \geq 5$$

\section{Leibniz Superalgebras with two-dimensional even part}
%\section{Case $n=1$ , $\{X_0,X_1\}+\{Y_1,Y_2,\dots,Y_m\}$}

According to the precedent cases, now the first subcase to
consider is $m=4$. In these cases, that is $n=2$ and $m\geq 4$, in
order to prove the following proposition we use the condition of
nilpotency, graded Leibniz identity and generic changes of basis
(isomorphisms). However, in this case  the complexity is higher
than the other cases.

\begin{prop} Let $L$ be
any Leibniz superalgebra $L \in {\cal ZF}^{2,m}$, with $m\geq 4$.
Then it is isomorphic to one of the following Leibniz
superalgebras, pairwise non-isomorphic, that can be expressed in
an adapted basis $\{X_0,X_1,Y_1,Y_2,Y_3,Y_4,\dots,Y_m\}$ by

$$\mu_{k} =\left\{\begin{array}{ll}
[X_0,X_0]=X_{1}& \\[1mm]
[Y_i,X_0]=Y_{i+1} & 1\leq i \leq m-1\\[1mm]
%[X_0,Y_1]=-Y_2 &  \\[1mm]
[X_0,Y_{j}]=-Y_{j+1}& 1 \leq j \leq m-1\\[1mm]
[Y_i,Y_j]=(-1)^{i+1} X_1& 1 \leq i,j \leq m,\
i+j=2k+2\\
\end{array}\right.$$
with $1\leq k \leq \lfloor\frac{m-1}{2}\rfloor$. $$\mu_{k}
=\left\{\begin{array}{ll}
[X_0,X_0]=X_{1}& \\[1mm]
[Y_i,X_0]=Y_{i+1} & 1\leq i \leq m-1\\[1mm]
[X_0,Y_1]=-Y_2 + Y_m&  \\[1mm]
[X_0,Y_{j}]=-Y_{j+1}& 2 \leq j \leq m-1\\[1mm]
[Y_i,Y_j]=(-1)^{i+1} X_1& 1 \leq i,j \leq m,\
i+j=2k+2-2 \lfloor\frac{m-1}{2}\rfloor\\
\end{array}\right.$$
with $\lfloor\frac{m}{2}\rfloor\leq k \leq m-2$.

$$\begin{array}{ll} \mu_{m-1} =\left\{\begin{array}{ll}
[X_0,X_0]=X_{1}& \\[1mm]
[Y_i,X_0]=Y_{i+1} & 1\leq i \leq m-1\\[1mm]
%[X_0,Y_1]=\a Y_m& \a \in \{0,1\} \\[1mm]
[Y_1,Y_1]=X_1 & \\[1mm]
\end{array}\right.&
\mu_{m} =\left\{\begin{array}{ll}
[X_0,X_0]=X_{1}& \\[1mm]
[Y_i,X_0]=Y_{i+1} & 1\leq i \leq m-1\\[1mm]
[X_0,Y_1]=Y_m&  \\[1mm]
[Y_1,Y_1]=X_1 & \\[1mm]
\end{array}\right.\\ \\
\mu_{m+1}=\left\{\begin{array}{ll}
[X_0,X_0]=X_{1}& \\[1mm]
[Y_i,X_0]=Y_{i+1} & 1\leq i \leq m-1\\[1mm]
%[X_0,Y_1]=-Y_2+ \a Y_m& \a \in \{0,1\}\\[1mm]
[X_0,Y_{j}]=-Y_{j+1} & 2\leq j \leq m-1\\[1mm]
[Y_1,Y_1]=X_1 & \\[1mm]
\end{array}\right.&
\mu_{m+2}=\left\{\begin{array}{ll}
[X_0,X_0]=X_{1}& \\[1mm]
[Y_i,X_0]=Y_{i+1} & 1\leq i \leq m-1\\[1mm]
[X_0,Y_1]=-Y_2+ Y_m& \\[1mm]
[X_0,Y_{j}]=-Y_{j+1} & 2\leq j \leq m-1\\[1mm]
[Y_1,Y_1]=X_1 & \\[1mm]
\end{array}\right.
\end{array}$$
\end{prop}

\begin{thm}  When $m$ is odd $(m\geq 4),$ then
$O(\mu_{\frac{m-1}{2}})$ is an open subset of ${\cal ZF}^{2,m}.$
\end{thm}

\begin{proof} The present proof is evident, because all the above Leibniz
superalgebras for $m$ even have nilindex $m,$ but for $m$ odd,
$\mu_{\frac{m-1}{2}}$ will have maximal nilindex $m+1$.
\end{proof}

\begin{lem}Let $L=L_0\oplus L_1$ be an arbitrary Leibniz filiform superalgebra
with $gz(L)=(1,1|4)$. Then $L_0$ is the abelian algebra and $L$
will have nilindex $4$.
\end{lem}

\begin{proof} The above lemma is a consequence of lemma \ref{lema}, of the
graded Leibniz identity and of the condition of being nilpotent.
\end{proof}

\begin{prop}$$f(2,4)=5$$
\end{prop}
\begin{proof} Consider $\mu_6 \oplus <Y_4>$, where $\mu_6$ is the
zero-filiform Leibniz superalgebra of maximal nilindex for the
case $n=2$ and $m=3$ (see Lemma \ref{lema2}). \end{proof}

We have as corollaries the following proposition and theorem.

\

\begin{prop} $${\cal M}^{2,4}\not \subset{\cal ZF}^{2,4}$$

\end{prop}

\

\begin{thm}
$$f(2,m)=m+1, \mbox{ if } m \mbox{ is odd and }m\geq 5$$
$$f(2,m)=m, \mbox{ if } m \mbox{ is even and }m\geq 6$$
\end{thm}

\

\begin{thm}  If $m \geq 5$, then
$${\cal M}^{2,m}\not \subset{\cal ZF}^{2,m}$$
\end{thm}

\begin{proof} For $m$ odd we can consider the Lie superalgebra $K^{2,m}$
(see \cite{GKN}) that has nilindex $m+1$ and it is not included in
${\cal ZF}^{2,m}$. For $m$ even we can consider $K^{2,m-1} \oplus
\CC$ with nilindex $m$. \end{proof}

\section{Conjecture}

\noindent{\bf Conjecture 2.}

If $n+1$ and $n$ are the dimensions of the even part and the odd
part, respectively, the only non split Leibniz superalgebra of
type $(n+1,n)$ and nilindex $2n$ is the following filiform Leibniz
superalgebra. The law of this superalgebra can be expressed, in an
adapted basis $\{X_0,X_1,\dots,X_n, Y_1,Y_2,\dots,Y_n\}$, by

%The only one Leibniz superalgebra, up
%to isomorphism, non split, of type $(n+1,n)$ ($n+1$, $n$ are the
%dimensions of the even part and odd part respectively) with
%nilindex $2n$ is the following filiform Leibniz superalgebra that
%can be expressed in an adapted basis
%$\{X_0,X_1,\dots,X_n\}+\{Y_1,Y_2,\dots,Y_n\}$ as

$$R^{n+1,n}=\left\{\begin{array}{ll}
[X_i,X_0]=X_{i+1}& 1 \leq i \leq n-1\\[1mm]
[X_0,X_0]=X_2 & \\[1mm]
[X_0,Y_1]=\frac{1}{2} Y_2 &\\[1mm]
[X_i,Y_1]=\frac{1}{2} Y_{i+1} & 1\leq i \leq n-1\\[1mm]
[Y_j,X_0]=Y_{j+1} & 1 \leq j \leq n-1\\[1mm]
[Y_1,Y_1]=X_0  &\\[1mm]
[Y_i,Y_1]=X_{i}& 1 \leq i \leq n \\[1mm]
\end{array}\right.$$

\bibliographystyle{amsplain}

\begin{thebibliography}{99}

\bibitem{AO0}
{\rm Albeverio S., Ayupov Sh.A., Omirov B.A.}{\it On nilpotent and
simple Leibniz algebras.} Rheinische
Friedrich-Wilhelms-Universitat Bonn, Preprint, 14 p, 2002
(published in Comm. in Algebra, v. 33(1), 2005, p. 159-172).

\bibitem{AOR}
{\rm Albeverio S., Omirov B.A., Rakhimov I.S.} {\it Varieties of
nilpotent complex Leibniz algebras of dimension less than five.}
Comm. in Algebra, v. 33(5), 2005, p. 1575-1585.

\bibitem{AO1}
{\rm Ayupov Sh.A., Omirov B.A.}  {\it On some classes of nilpotent
Leibniz algebras.} Siberian Math. Journal, v. 42(1), 2001, p.
18--29.

\bibitem{AO2}
{\rm Ayupov Sh.A., Omirov B.A.} {\it On Leibniz algebras.} Algebra
and operators theory, Proceedings of the Colloquium in Tashkent
1997. Kluwer Academic Publishers, 1998, p. 1--13.

\bibitem{BS}
{\rm Burde D., Steinhoff C.}  {\it Classification of orbit
closures of 4-dimensional complex Lie algebras.} J. of Algebra, v.
214, 1999, p. 729--739.

\bibitem{Jac}
{\rm Jacobson N.} {\it Lie algebras.} Interscience Publishers,
Wiley, New York, 1962.

\bibitem{GO}
{\rm Grunewald F., O'Halloran J.} {\it Varieties of nilpotent Lie
algebras of dimension less than six.} J. of Algebra, v. 112(2),
1988, p.315--326.

\bibitem{Liv}
{\rm Livernet M.}  {\it Rational homotopy of Leibniz algebras.}
Manuscripta math., 96, 1998, p. 295--315.

\bibitem{Loday}
{\rm Loday J.-L.} {\it Une version non commutative des
alg$\acute{e}$bres de Lie: les alg$\acute{e}$bres de Leibniz.}
Ens. Math., 39, 1993, pp. 269--293.

%\bibitem{Gilg}
%{\rm Gilg M.} {\it Super-alg$\acute{e}$bres.}
%\par PhD thesis, University of Haute Alsace, 2000.

\bibitem{GKN}
{\rm G\'{o}mez J. R., Khakimdjanov Yu., Navarro R.M.} {\it Some
problems concerning to nilpotent Lie superalgebras.} J. Geometry
and Physics, v. 51(4), 2004, p. 472--485.


\bibitem{GK}
{\rm Goze M., Khakimdjanov Yu.}  {\it Nilpotent Lie Algebras.}
 Kluwer Academic publishers, Dordrecht, 1996, 336 p.

%\bibitem{Nav}
%{\rm Navarro R.M.} {\it Superalgebras de Lie nilpotentes.}
%\par PhD thesis, University of Sevilla, 2001.

\end{thebibliography}

\end{document}